\documentclass[11pt,a4paper, DIV=12]{scrartcl}
\usepackage[utf8]{inputenc}
\usepackage[T1]{fontenc}
\usepackage{amsmath}
\usepackage{amssymb}
\usepackage{amscd}
\usepackage{amsthm}
\usepackage{proof}
\usepackage{color}
\usepackage{enumerate}
\usepackage{hyperref}
\setkomafont{sectioning}{\bfseries}

\newtheorem{theorem}{Theorem}[section]
\newtheorem{coro}[theorem]{Corollary}
\newtheorem{lemma}[theorem]{Lemma}
\newtheorem{prop}[theorem]{Proposition}

\newtheorem{definition}[theorem]{Definition}

\theoremstyle{definition}
\newtheorem*{remark}{Remark}

\newcommand{\exend}{\hfill $\Diamond$}

\newcommand{\ts}{\hspace{0.5pt}}

\newcommand{\nfunction}{\mathcal{N}}

\newcommand{\NN}{\mathbb{N}}
\newcommand{\ZZ}{\mathbb{Z}}
\newcommand{\CC}{\mathbb{C}}
\newcommand{\RR}{\mathbb{R}}
\newcommand{\TT}{\mathbb{T}}

\title{Spectral theory of dynamical systems as diffraction theory of
sampling functions}

\author{D.~Lenz\footnote{ Mathematisches Institut, Friedrich Schiller
Universit{\"a}t Jena, 07743 Jena, Germany, daniel.lenz@uni-jena.de},
}

\begin{document}

\maketitle

\begin{abstract}
We consider topological dynamical systems over  $\ZZ$ and, more
generally, locally compact, $\sigma$-compact abelian groups. We
relate spectral theory and diffraction theory.  We first use a a
recently developed general framework of diffraction theory to
associate an autocorrelation and a diffraction measure to any
$L^2$-function over such a dynamical system. This diffraction
measure is shown to be the spectral measure of the function. If the
group has a countable basis of the topology one can also exhibit the
underlying autocorrelation by sampling along the orbits. Building on
these considerations we then  show how the spectral theory of
dynamical systems can be reformulated via diffraction theory of
function dynamical systems. In particular, we show that the
diffraction measures of suitable factors provide a complete spectral
invariant.
\end{abstract}


\section*{Introduction}
During the last three decades,  the mathematical theory of aperiodic
order has become a field of substantial interest.  This is not the
least due to the discovery (honored with a Noble Prize in Chemistry)
of certain materials  - later known as quasicrystals - exhibiting
this form of order \cite{Danny}. The discovery of such substances
came as a complete surprise to physicists and materials scientists
alike. These materials are characterized by their remarkable
diffraction properties: They exhibit pure point diffraction
(indicating long range order) and at the same time their diffraction
patterns exhibit symmetries which are incompatible with a lattice
structure. Hence, their structure exhibits   long range aperiodic
order.

The mathematical  study of diffraction of quasicrystals and,  more
generally, the mathematical study of aperiodic order has profited
tremendously from ideas and methods from dynamical systems and
stochastic processes, see e.g. the survey collections
\cite{BaakeMoody00,KLS,TAO2} and the monograph \cite{TAO}. On the
other hand, recent years have also seen a flow of ideas in the other
direction, i.e. from the study of aperiodic order to dynamical
systems.
 In fact, inspired by diffraction theory for quasicrystals,  a
characterization for pure discrete spectrum is given for general
dynamical systems in \cite{Len}. Similarly, the proof of recent
results of regularity of  tame systems was paved by considerations
on certain models of aperiodic order \cite{FGJO}.

Another instance of a flow of ideas from the study of  aperiodic
order to the treatment of general dynamical systems concerns
spectral theory of dynamical systems. This is the topic of the
present article.  Specifically, it is possible to consider the
spectral measures arising in the study of dynamical systems as
diffraction measures. Certainly, this will not come as a surprise to
the experts. In fact, for  special systems an indication of this is
already given in \cite{BLE}. However, it does not seem to be
discussed explicitly and in any form of  detail for general
dynamical systems in the existing literature. As this point of view
may be useful for further developments,  we present a rather
complete discussion here. In this discussion  the connection between
spectral theory and diffraction theory is investigated along two
different (but related) lines.

The first line is provided in the first part of the article starting
with Section \ref{sec-framework}. There, we  show that the framework
for diffraction theory developed in \cite{LM} is general enough to
not only cover the `usual' situations treated in diffraction theory
but also allows one to recover the spectral measures as diffraction
measures. As a by-product,  we obtain in a rather simple way a
structural understanding of (generalizations of) recent results of
\cite{KMSS}. These considerations  can be cast in the framework of
dynamical systems over general locally compact, $\sigma$-compact
abelian groups. To ease the presentation and as this is a
particularly relevant case,   we have decided to  first present the
case of dynamical systems over the group of integers  in Section
\ref{sec-framework}.  Subsequently,  we consider the case of general
locally compact, $\sigma$-compact abelian groups in Section
\ref{sec-general}.  Most of the considerations  concerning the
integers can also be adapted to treat actions of the semigroup of
natural numbers rather than the integers, see also \cite{Q} for a
treatment of related problems on one-sided subshifts. We leave the
details to the interested  reader.

Diffraction theory was originally developed for dynamical systems of
point sets and generalizations thereof. In this context, the most
general framework is provided by the  diffraction theory for measure
dynamical systems developed in \cite{BL,LS}. The considerations
discussed above  show that on the formal level diffraction theory
for measure dynamical systems and spectral theory of general
dynamical systems can be treated on the same footing. This also begs
the question whether there is an intrinsic connection between
diffraction theory for measure dynamical systems and spectral theory
for general dynamical systems. An affirmative answer to this
question constitutes the  second line of  connecting diffraction
theory and spectral theory. The corresponding discussion is given in
the second part of this article. This starts in Section
\ref{sec-comparison} and  builds up on the first part of the
article. In  Section \ref{sec-comparison}  we recall basic the
diffraction theory for measure dynamical systems following
\cite{BL}. Moreover, we introduce a special class of measure
dynamical systems, termed function dynamical systems, which are of
particular relevance for our subsequent considerations. We then show
in Section \ref{sec-factors} that any dynamical system has canonical
factors which are function dynamical systems and that the
diffraction of these function dynamical systems encodes the spectral
theory of the original system. On the structural level this is our
main result. It may be summarized as follows (compare Theorem
\ref{t:spectralmeasure-as-diffraction-of-TMDS} and Proposition
\ref{prop:diffraction-factor} for precise statements):

\medskip

\textbf{Result.}  \textit{The spectral theory of a general dynamical
systems is   the  diffraction theory of those of its  factors which
are function dynamical systems}.

\medskip

We can use this   to show for any dynamical system that  the
diffraction measures of its function dynamical system factors form a
complete spectral invariant (Theorem \ref{t:analogue-BLE}) and so do
the diffraction measures of its translation bounded measure
dynamical system factors  (Corollary \ref{cor:analogue-BLE}). This
type of result has so far only be known for a very restricted class
of measure dynamical systems, i.e. Delone dynamical systems of
finite local complexity \cite{BLE}. Beyond the structural insights,
our results may also be of interest for specific questions.  In
fact, they allow us to study arbitrary dynamical systems via methods
developed for the study of measure dynamical systems. As a concrete
application we provide in Section \ref{s:application}  a  criterion
for discrete spectrum based on recent investigations of measure
dynamical systems carried out in \cite{SpStr}. This can be seen as a
variant of a well-known criterion for $\ZZ$.

We also point out that our results give that  any spectral measure
(of a continuous function) is a diffraction measure. Hence, they
allow one to generate - in a perfectly natural way - interesting
examples of diffraction measures.

The preceding results are set in the topological category in that
they  deal with continuous functions and factor maps. It is also
possible to deal with their analogues in a measurable setting.  This
is discussed in the Section \ref{sec-measurable}.

Finally, for the convenience of the reader we include in an appendix
a short review of basics of the diffraction theory developed in
\cite{LM}.

\bigskip

\textbf{Acknowledgements.} The author would like to thank Michael
Baake for helpful and encouraging  comments on the manuscript and
Nicolae Strungaru for enlightening discussions on abstract versions
of diffraction theory and mean almost periodicity. The author is
also indebted to  inspiring discussions with Gerhard Keller,
Christoph Richard and Jeong-Yup Lee at the workshop on 'Model Sets
and Aperiodic Order' in Durham (2018) and he would like to express
his thanks to the organizers of this workshop. The  manuscript was
completed during a stay at the university of Geneva and hospitality
of the department of mathematics is gratefully acknowledged.

\section{Dynamical systems over $\ZZ$}\label{sec-framework}
In this section,  we first introduce our main actors. These are
dynamical systems and their spectral theory as well as diffraction
theory. Throughout the article all dynamical systems are topological
dynamical systems.

Diffraction theory has been developed  in various contexts in
various degrees of generality  in the last decades, see e.g. the
survey \cite{BL} for a gentle introduction. Here, we follow the most
general framework given in \cite{LM} (see Appendix as well). In
fact, the framework of \cite{LM} can be slightly simplified for
actions of the group of integers and this is the setting we present.
One of these simplifications is that we can deal with functions
throughout and do not have to consider measures. Having conveniently
set up the framework,   the derivation of the main results is then
straightforward.

\bigskip

Throughout we denote the group of integers by $\ZZ$ and  the vector
space of functions on $\ZZ$ with finite support by $C_c (\ZZ)$. For
any $n\in\ZZ$ the  characteristic function of $\{n\}\subset \ZZ$ is
denoted by $1_n$. The convolution between $h : \ZZ\longrightarrow
\CC$ and $\varphi \in C_c (\ZZ)$ is the function
$$ h \ast \varphi : \ZZ\longrightarrow \CC, n\mapsto \sum_{k\in\ZZ}
h (n-k) \varphi (k).$$ For $\varphi \in C_c (\ZZ)$ we define
$\widetilde{\varphi} \in C_c (\ZZ)$ via
$$\widetilde{\varphi} (n) =
\overline{\varphi (-n)}$$ for all $n\in \ZZ$. A  function $p$ on
$\ZZ$ is \textit{positive definite} if
$$\sum_{n\in \ZZ} p (n) (\varphi \ast \widetilde{\varphi})(n) \geq
0$$ for all $\varphi \in C_c (\ZZ)$. By a result of Bochner any
positive definite function $p$ is the Fourier transform of   a
unique measure $\varrho$ on the unit circle
$$\TT:=\{z\in \CC: |z| =1\}, $$
i.e.
$$p(n) =\int_{\TT} z^n d\varrho (z)$$
for all $n\in\ZZ$, see e.g. \cite{BF}.

By a  \textit{dynamical system} we mean a triple $(X,\alpha, m)$
consisting of a compact space $X$,  a continuous action
$$\alpha : \ZZ\times X \longrightarrow X$$
of $\ZZ$ on $X$ and an $\alpha$-invariant probability measure $m$ on
$X$. Clearly, this action is completely determined by $\alpha_1$ as
we have $\alpha_n = (\alpha_1)^n$ for any $n\in\ZZ$. The dynamical
system is \textit{ergodic} if any $\alpha$-invariant measurable
subset of $X$ has measure zero or one.

Any dynamical system $(X,\alpha, m)$ comes naturally with the
Hilbert space $L^2 (X,m)$ with inner product
$$\langle f, g\rangle := \int_X f \overline{g}  d\mu$$
as well as a unitary map $U$ on  $L^2 (X,m)$ defined via
$$U f = f \circ \alpha_{-1}.$$
The unitary $U$ is known as \textit{Koopman operator}. Such
operators play an enormous  role in the study of dynamical systems,
see e.g. the recent monograph \cite{EFHN} for the (ever increasing)
interest in this type of operators.

Consider now  an  $f\in L^2 (X,m)$. Then, a short computation shows
$$\sum_{n\in\ZZ} \langle U^n f, f\rangle (\varphi \ast
\widetilde{\varphi}) (n) = \|\sum_n \varphi (n) U^n f\|^2 \geq 0$$
for all $\varphi \in C_c (\ZZ)$.  Hence,
 the function  $$
\ZZ\longrightarrow \CC, n\mapsto \langle U^n f, f \rangle,$$ is
positive definite. By the theorem of Bochner mentioned above, there
exists  then a unique measure $\varrho^f$ on the unit circle with
$$\langle U^n f, f\rangle = \int_{\TT} z^n d\varrho^f (z)$$
for all $n\in \ZZ$.   This measure is called the \textit{spectral
measure} of $f$, see e.g. \cite{Ped}. Spectral theory of dynamical
systems is the study of these spectral measures.

We now turn to diffraction theory. It  is developed to describe the
outcome of diffraction experiments. Here, the piece of matter to be
analyzed in the diffraction experiment is   modeled by a function
and the  diffraction measure of this function describes the outcome
of the diffraction experiment. A discussion of the physics behind
this can be found in \cite{Cow}. The mathematical side is developed
in the fundamental work \cite{Hof}, see  the monograph \cite{TAO} as
well for a general discussion. Specifically, the
\textit{autocorrelation} of a function $h : \ZZ\longrightarrow \CC$
is defined as the pointwise limit of the functions
$$\frac{1}{2 n} \sum_{- n \leq k,l \leq n } h(k) \overline{ h (l)} 1_{k-l}$$ if
this limit exists. This autocorrelation can easily be seen to be
positive definite. Hence, by Bochner theorem,  it is the Fourier
transform of a positive measure  on $\TT$.  This measure is called
the \textit{diffraction measure} of $h$.

\begin{remark} Diffraction theory in one dimension is often
developed for measures on $\RR$. In this section, we rather deal
with functions on $\ZZ$.  It is possible to consider functions on
$\ZZ$ as measures on $\RR$ by associating to $h : \ZZ\longrightarrow
\RR$ the measure $\delta_h:=\sum_{x\in \ZZ} h(x) \delta_x$ (with
$\delta_p$ being the unit point mass at $p$). The autocorrelation
and diffraction measure of $h$ and of $\delta_h$ are related. The
diffraction measure of $\delta_h$  is a periodic extension (with
period $1)$ of the diffraction measure on $h$, where $\TT$ is
considered as $[0,1)\subset \RR$ , see \cite{BL,BLE} for details.
\exend
\end{remark}

In typical diffraction  situations,  one is not given one function
but rather a family of functions arising from a dynamical system.
This is captured by the somewhat more involved framework of
diffraction arising  from dynamical systems. In order to set this up
one needs a dynamical system together with one more ingredient. Let
$\beta$ be the natural action of $\ZZ$ on $C_c (\ZZ)$ i.e. $(\beta_n
\varphi) (k) = \varphi (k - n)$. The mentioned ingredient is then  a
map
$$\nfunction : C_c (\ZZ) \longrightarrow L^2(X,m)$$
satisfying the following two properties (see Appendix as well):

\begin{itemize}

\item[(N1)] $\nfunction$ is linear.

\item[(N2)] $\nfunction$ is equivariant (i.e.
$\nfunction(\beta_n \varphi) = U^n \nfunction(\varphi)$ for all
$n\in\ZZ$ and $\varphi \in C_c (\ZZ)$).

\end{itemize}

Any such map then comes with a unique function  $\gamma =
\gamma^{(\nfunction)}$ on $\ZZ$ satisfying
$$\sum_{n\in\ZZ} \gamma (n) (\varphi \ast \widetilde{\psi}) (n)   = \langle \nfunction(\varphi),
\nfunction(\psi)\rangle$$ for all  $\varphi,\psi \in C_c (\ZZ)$.
Indeed, a short calculation confirms that $\gamma$ is given by
$\gamma = \sum_{n} c_n 1_n$ with
$$c_n := c_n^{(\nfunction)}:=\langle \nfunction (1_n), \nfunction(1_0)\rangle.$$
Note that $c_n$ satisfies
$$c_n = \langle U^n  \nfunction (1_0),  \nfunction (1_0)\rangle$$
by equivariance of $\nfunction$. The function  $\gamma$ is called
the \textit{autocorrelation} of $\nfunction$. Again, $\gamma$ is
easily be seen to be  positive definite. Hence, it has a Fourier
transform, which is a measure on $\TT$. This measure is called the
\textit{diffraction measure} of $\nfunction$.

This immediately gives the following.

\begin{lemma}[Autocorrelation and spectral measure]\label{lemma-main} Let $(X,\alpha,m)$ be a dynamical system and
$\nfunction : C_c (\ZZ)\longrightarrow L^2 (X,m)$ satisfy (N1) and
(N2). Then, the spectral measure of $N(1_0)$ is the diffraction
measure of $\nfunction$.
\end{lemma}
\begin{proof}  By the preceding discussion we have
$$\gamma^\nfunction (n) = c_n = \langle U^n \nfunction
(1_0),\nfunction (1_0)\rangle = \int z^n
d\varrho^{\nfunction(1_0)}$$ for all $n\in\ZZ$. Now, the result
follows by taking Fourier transforms on both sides.
\end{proof}


\begin{remark} In the general situation treated in \cite{LM} existence of $\gamma$
is not automatically satisfied. For $\ZZ$  (and more generally for
discrete groups) it is automatically satisfied.\exend
\end{remark}

From the defining properties of $\nfunction$  it is not hard to see
that $\nfunction(1_0)$ is an $L^2$-function and the map $\nfunction$
is determined by $\nfunction(1_0)$. In fact - as will be discussed
more thoroughly below in this section -  any  element from $ L^2
(X,m)$ gives rise to such a map $\nfunction$. So, such functions
$\nfunction$ are  in one-to-one correspondence with functions in
$L^2 (X,m)$.

Under a suitable ergodicity assumption,  it is possible to
understand the  autocorrelation of $\nfunction$ as the
autocorrelation of the sampling along $\nfunction(1_0)$, i.e. as the
autocorrelation of the \textit{sampling functions}
$$\ZZ\longrightarrow \CC, n\mapsto
\nfunction(1_0) (\alpha_{-n} (x)),$$ for $x\in X$.
 This is discussed next.

\begin{lemma}\label{sampling} Let $(X,\alpha,m)$ be ergodic and $\nfunction : C_c
(\ZZ)\longrightarrow L^2 (X,m)$ as above. Set $f:= \nfunction(1_0)$.
 Then, for $m$-almost every $x\in X$, the autocorrelation
$\gamma^\nfunction$ is equal to the autocorrelation of the function
$\ZZ\longrightarrow \CC, n\mapsto f(\alpha_{-n} (x))$.
\end{lemma}
\begin{proof} This follows by a standard computation from Birkhoff's ergodic theorem. We include it
for the convenience of the reader: We have to show that the
functions  $\gamma_n^x$ with
$$\gamma_n^x  = \frac{1}{2 n} \sum_{- n \leq j,l \leq n }
f(\alpha_{-j} (x))\overline{f(\alpha_{-l} (x))} 1_{j-l}$$ converge
pointwise to $\gamma^{(\nfunction)}$ for almost every $x\in X$. Fix
an index $k$. Then, the coefficient of $1_k$ in
$\gamma^{(\nfunction)}$ is given by
$$\langle \nfunction (1_k), \nfunction (1_0)\rangle = \langle f \circ \alpha_{-k}, f\rangle.$$
On the other hand, the coefficient of $1_k$
 in
$\gamma_n^x$ for sufficiently large $n$ is essentially given by
$$\frac{1}{2n} \sum_{-n \leq l \leq n} f(\alpha_{-k-l} (x))
\overline{f(\alpha_{-l} (x))}.$$ By Birkhoff's ergodic theorem, this
can easily be seen to converge for $m$-almost every $x\in X$ to
$$\int_X f(\alpha_{-k} x) \overline{f ( x)} d m  (x) =  \langle f \circ \alpha_{-k}, f \rangle.$$
As there are only countably many $k\in\ZZ$, we find pointwise
convergence of the functions in question for $m$-almost every $x\in
X$. This finishes the proof.
\end{proof}

\begin{remark} (a) Note that the previous lemma provides a connection
between the two versions of diffraction theory discussed above, viz
the version based on individual functions and the version based on
dynamical systems.

(b) If the system is uniquely ergodic (i.e. there exists only one
$\alpha$-invariant probability measure on $X$) and $\nfunction(1_0)$
is Riemann integrable, the limit in the preceding lemma can easily
be shown to exist for all $x\in X$ by Oxtoby's theorem. \exend
\end{remark}

\begin{remark}[Situation for real $f$] For real valued  $f$ the map
$\nfunction$ is also real (i.e. satisfies $\nfunction
(\overline{\varphi}) =\overline{\nfunction(\varphi)}$ for all
$\varphi \in C_c (\ZZ)$).  This brings an additional symmetry into
play. More specifically, it implies  that  $\gamma$ is real-valued.
Now,  $\gamma$ is positive definite and, hence, satisfies $\gamma
(n) = \overline{\gamma (-n)}$ for all $n\in\ZZ$. Taken together
this then gives that $\gamma$ is symmetric under reflection $R :
\ZZ\longrightarrow \ZZ,n \mapsto -n$. In particular, $\gamma$ is
also the autocorrelation of the dynamical system $(X,\alpha^R,m)$
with respect to $\nfunction^R$, where $\alpha^R$ is the reflected
action defined by $\alpha^R_n = \alpha_{-n}$ and $\nfunction^R$ is
defined by $\nfunction^R (\varphi) :=\nfunction (\varphi \circ R)$.
In particular, $\gamma^\nfunction$ also arises by sampling along the
function $n\mapsto f(\alpha_n (x))$ for almost every $x\in X$.
\exend
\end{remark}

Let now  $(X,\alpha,m)$ be a dynamical system. Let $f\in L^2(X,m)$
be given and   define $\nfunction^f : C_c (\ZZ)\longrightarrow L^2
(X,m)$ via
$$\nfunction^f (\varphi) := \sum_{n\in\ZZ} \varphi (n) U^n f.$$
It is not hard to see that $\nfunction^f$ satisfies the properties
(N1) and (N2).  Hence, there is an autocorrelation and a diffraction
associated to $\nfunction^f$. We denote this autocorrelation by
$\gamma^f$ and this diffraction by $\omega^f$ and call it the
\textit{autocorrelation} and \textit{diffraction measure} of $f$,
respectively. Clearly,
$$\nfunction^f (1_0) = f$$
holds.  Hence, from Lemma \ref{lemma-main} we immediately derive the
main result of this section.

\begin{theorem}[Spectral measure as diffraction measure]
 Let $(X,\alpha,m)$ be a dynamical system and $f\in
L^2 (X,m)$ be given. Then, $\omega^f = \varrho^f$.
\end{theorem}

If the system is ergodic,  it is  possible to compute the
autocorrelation almost surely as the  autocorrelation of the
sequences arising from sampling along $f$. Indeed, we have the
following immediate consequence of the previous theorem and Lemma
\ref{sampling}.

\begin{coro}[Pointwise sampling] Let $(X,\alpha,m)$ be ergodic and $f\in L^2 (X,m)$
be real valued. Then, for $m$-almost every $x\in X$, the
autocorrelation $\gamma^f$ is equal to the autocorrelation of the
function  $\ZZ\longrightarrow \RR, n\mapsto f(\alpha_{-n} (x))$.
\end{coro}

\begin{remark} If the dynamical system is uniquely ergodic, then the
autocorrelation exists for all $x\in X$, whenever $f$ is continuous
(or, more generally, Riemann integrable).\exend
\end{remark}

\begin{remark}[Relationship with \cite{KMSS}]
The previous result immediately implies that dynamical systems whose
spectral measures have certain properties  will give rise to
diffraction measures with these properties  as well.

In particular, any dynamical system with pure discrete spectrum,
i.e. pure point spectral measures, will give rise to pure point
diffraction measures $\omega^f$ associated to its real valued
functions $f$. Such systems are at the core of the theory of
aperiodic order (see e.g. \cite{TAO,BL,KLS} for recent introductory
and survey type treatments  of this topic).

Likewise, a dynamical system with sufficiently strong mixing
properties to ensure absolute continuity of its spectral measures
will have diffraction measures $\omega^f$ which are absolutely
continuous (up to a possible atom at $0$). Such systems include
Bernoulli systems and have also attracted attention in the context
of diffraction theory, see e.g. last chapter of \cite{TAO} or the
survey article \cite{BBM}.

This provides also a structural understanding of the phenomena
arising in the examples given in \cite{KMSS}. Indeed, the examples
treated there are exactly dynamical systems with pure point spectrum
and dynamical systems with strong mixing properties respectively.
Note, however, that \cite{KMSS} deals also with non-invertible
systems, whereas we have restricted our attention  to invertible
systems. Finally, let us point out that existence of sampling along
orbits is also already discussed in \cite{KMSS}. \exend
\end{remark}

\section{Dynamical systems over general locally compact abelian groups}\label{sec-general}
In the last section,  we have considered dynamical systems with
actions of the group of  integers. Here, we consider dynamical
systems over general locally compact, $\sigma$-compact, abelian
groups and indicate how the theory developed above carries over to
this case as well. This requires some further technical elaborations
not present in the case of $\ZZ$.

\bigskip

We consider  a compact space $X$ equipped with a continuous  action
$$G\times X\longrightarrow X, (t,x)\mapsto
\alpha_t (x),$$ of a locally compact, $\sigma$-compact,  abelian
group $G$ and a probability measure $m$, which is invariant under
the action of $G$. We then call $(X,G,m)$ a \textit{dynamical system
over the space $X$} and write $tx$ instead of $\alpha_t (x)$ for
$t\in G$ and $x\in X$. Likewise the composition on $G$ itself is
written additively.

The dual group of $G$ (i.e.the set of all continuous group
homomorphisms from $G$ to the unit circle $\TT$) is denoted as
$\widehat{G}$. Haar measure on $G$ is denoted as $m_G$. The  space
of continuous functions on $G$ with compact support is denoted by
$C_c (G)$ and the convolution between a measure $\eta$ on $G$ and a
$\varphi \in C_c (G)$ is the function $\eta \ast \varphi :
G\longrightarrow \CC$ with
$$\eta \ast \varphi (t) = \int_G \varphi (t - s) d\eta (s)$$
for $t\in G$. A function $F$ on $G$  (whose restrictions to compact
sets belong to $L^1 (G,m_G)$) may be considered as the measure $F
m_G$ and we can accordingly define $F\ast \varphi := (F m_G) \ast
\varphi$ for $\varphi \in C_c (G)$. For $\varphi \in C_c (G)$ we
define $\widetilde{\varphi} \in C_c (G)$ by $\widetilde{\varphi}(t)
= \overline{\varphi(-t)}$ for $t\in G$.

The action of $G$ on $X$ induces  unitary operators $T_t : L^2
(X,m)\longrightarrow L^2 (X,m)$ with $$ T_t f  = f((-t)\;\cdot)$$
for each $t\in G$. Here, $L^2 (X,m)$ is the Hilbert space of
(equivalence classes of) square integrable functions on $X$. It is
equipped with the inner product
$$\langle f, g\rangle = \int_X  f \overline{g} \,  dm$$
and the associated norm
$$\|f\| := \|f\|_2:=\sqrt{\langle f, f\rangle}$$
for $f,g\in L^2 (X,m)$. Then, for  each $f\in L^2 (X,m)$  there
exists a unique finite measure $\varrho^f$ on $\widehat{G}$ with
$$\langle T_t f,  f\rangle = \int_{\widehat{G}} \eta (t) d\varrho^f
(\eta)$$  for all $t\in G$. This measure is called the
\textit{spectral measure of $f$}.

Let now  $(X,G,m)$ be a dynamical system and consider  $f\in L^2
(X,m)$. We can then define the map
$$\nfunction^f : C_c (G)\longrightarrow L^2(X,m) \; \mbox{via} \;\:
\nfunction^f (\varphi):= \int \varphi (t) T_t f d m_G (t).$$ Note
that $\nfunction^f$ maps indeed in $L^2$ with  the bound
$$\|\nfunction^f (\varphi)\| \leq \int_G |\varphi(t)| \|T_t f\| dt =
\|f\| \|\varphi\|_{L^1 (G,m_G)}$$ for all $\varphi \in C_c (G)$.
 Clearly, this map  is linear,  and equivariant  and thus
satisfies the conditions (N1) and (N2)   given in the appendix. In
fact, it is not hard to see that it also satisfies (N3). This is
discussed next.

We have to show that  there exists a (unique) measure $\gamma^f$ on
$G$ with
$$\gamma^f (\varphi \ast \widetilde{\psi}) = \langle \nfunction^f (\varphi),
\nfunction^f (\psi)\rangle$$ for all  $\varphi,\psi\in C_c (G)$. In
fact, $\gamma^f$ is just given by  $(t\mapsto  \langle T_t f,
f\rangle) m_G$, i.e.  satisfies $$\int \varphi (t) d\gamma^f (t) =
\int \varphi(t) \langle T_t f, f\rangle dt$$ for all $\varphi \in
C_c (G)$. Indeed, unwinding the definitions and using the invariance
of Haar measure we find
\begin{eqnarray*} \int_G
(\varphi \ast \widetilde{\psi}) (t) \,  \langle T_t f, f\rangle dm_G
(t) &=&\int_G \left(\int_G \varphi (t+s) \overline{\psi (s)} dm_G
(s)\right) \langle T_{t} f, f\rangle \; dm_G(t)\\
&=&\int_G \int_G \varphi (t+s) \overline{\psi (s)}   \langle T_{t}
f, f\rangle \; dm_G(t) d m_G (s)\\
&=& \int_G \int_G \varphi (t) \overline{\psi (s)} \langle T_{t-s} f,
f\rangle \; dm_G(t) dm_G (s)\\
&=&\langle \nfunction^f (\varphi), \nfunction^f (\psi)\rangle
\end{eqnarray*}
for all $\varphi,\psi\in C_c (G)$. As this is just the defining
property of $\gamma^f$, we obtain the claimed respresentation.

From this formula for $\gamma^f$ we can also directly see  that it
is the Fourier transform of the spectral measure. More specifically,
for any $\varphi\in C_c (G)$ we obtain
\begin{eqnarray*}
\gamma^f (\varphi) &=&\int_G \varphi (t) \langle T_{t} f, f\rangle
\; dm_G(t)  \\
&=&\int_G \varphi (t)
\left(\int_{\widehat{G}} \eta(t) d\varrho^f (\eta)\right) dm_G (t)\\
&=& \int_{\widehat{G}} \left(\int_G \varphi (t) \eta (t) dm_G
(t)\right) d\varrho^f (\eta)\\
&=&
 \int_{\widehat{G}} \check{\varphi} \; d\varrho^f.
\end{eqnarray*}
Here, the inverse  Fourier transform $\check{\sigma}$ of $\sigma :
{G}\longrightarrow \CC$ is defined by
$$\check{\sigma} (\eta) = \int {\eta (t)} \sigma (t) d m_G(t).$$
As a consequence of these considerations we find in particular
$$\gamma^f (\varphi \ast \widetilde{\psi}) =\int_{\widehat{G}} {\check{\varphi}}
\; \overline{\check{\psi}} \;  d\varrho^{f}$$ for all
$\varphi,\psi\in C_c (G)$. Taking Fourier transforms on both sides
one obtains the following result (compare the argument in the last
section as well).

\begin{theorem}[Spectral measure as diffraction
measure]\label{thm-main-abstract}
 Let $(X,G,m)$ be a dynamical system and $f\in L^2
(X,m)$ be given. Then, the Fourier transform of $\gamma^f$ is the
spectral measure $\varrho^f$ of $f$.
\end{theorem}

\section{Diffraction theory of measure dynamical systems}
\label{sec-comparison} Diffraction theory is usually developed for
point sets or, more generally, measures. The leads then to
diffraction theory for measure dynamical systems, as developed in
\cite{BL}.  The first aim of this section is to briefly recall this
theory. A notable difference of this theory to the considerations
above is that there is no need to fix a function. Instead the
dynamical system itself already provides enough information to
define autocorrelation and diffraction.  This is quite remarkable
and the second aim of this section is to show how this can be
understood in light of the theory developed above via limiting
procedures. Moreover,  we present a relevant class of measure
dynamical systems for subsequent considerations. This  class -
called function dynamical systems below -  does not seem to have
been a systematic focus in the development of diffraction theory so
far.

\bigskip

Throughout this section $G$ is a locally compact, $\sigma$-compact
abelian group.  A measure $\mu$ on $G$ is called \textit{translation
bounded} if its total variation $|\mu|$ satisfies
$$\sup |\mu| (t + U) < \infty$$
for one (all) relatively compact open $U$ in $G$. We denote that set
of all  translation bounded measures  by $M^\infty (G)$ and equip it
with the vague topology. Then, $G$ admits a  natural action on
$M^\infty (G)$ by translations. More specifically, for $t\in G$ and
$\mu\in M^\infty (G)$ the measure  $t\mu $ is defined by
$$t\mu (\varphi) = \mu (\varphi (\cdot + t))$$
 for all $\varphi \in C_c (G)$.

Whenever  $X$ is a compact subset  of $M^\infty (G)$,  which is
invariant under the translation action to  and $m$ is an invariant
probability measure on $X$, we call $(X,G,m)$ a \textit{dynamical
system of translation bounded measures} or just TMDS for short. Such
a system comes with a canonical  map
$$\nfunction : C_c (G)\longrightarrow C (X), \varphi \mapsto \nfunction_\varphi  \mbox{ with } \nfunction_\varphi (\mu)
= \mu(\varphi) :=\int \varphi (s) d\mu (s).$$ Let us emphasize that
the existence of such a map is a distinctive feature of TMDS
compared to general dynamical systems.

It is not hard to see that $\nfunction$ is linear and equivariant
i.e. satisfies (N1) and (N2). Moreover,  as shown in \cite{BL}, it
also satisfies (N3), viz.  there exists then a unique translation
bounded measure $\gamma = \gamma^m$ on $(X,G,m)$ with
\begin{equation*}
\gamma (\varphi \ast \widetilde{\psi}) = \langle \nfunction_\varphi,
\nfunction_\psi\rangle
\end{equation*} for all $\varphi,\psi \in
C_c (G)$ and all $t\in G$.\footnote{Note that \cite{BL} uses a
different sign in the definition of $\nfunction$ (called $f$ there)
as well as has the inner product linear in the second argument. This
results in a different display of the formula for $\gamma$, viz.
$\gamma \ast \widetilde{\varphi}\ast \psi (0) = \langle f_\varphi,
f_\psi\rangle$. } In fact, \cite{BL} even contains an explicit
formula for $\gamma^m$, see \cite{Gou} as well). The measure
$\gamma$ is called the \textit{autocorrelation} of the TMDS. This
measure allows for a Fourier transform $\widehat{\gamma}$ which is a
(positive) measure on $\widehat{G}$. It is known as
\textit{diffraction} of the TMDS. By the previous discussion and
basic results on Fourier transform the diffraction satisfies
$$ |\check{\varphi}|^2 \widehat{\gamma} =
\varrho^{\nfunction_\varphi}$$ for all $\varphi \in C_c (G)$.

In the last sections, autocorrelation and diffraction arose from a
dynamical system together with a chosen function. Thus, it is
remarkable that for TMDS there is no need to chose a function.
Instead  the dynamical system alone  suffices as piece of data.
Intuitively speaking,  the reason is that any approximate unit in
$C_c (G)$ can be understood as providing a canonical choice of
`function'. More precisely, it is possible to exhibit the maps
$\nfunction$ and the Fourier transform $\widehat{\gamma}$ of the
autocorrelation measure $\gamma$ by limiting procedures which
involve chosen functions stemming from an approximate unit. Here are
the details:  If $(X,G,m)$ is a TMDS, we can consider an approximate
unit in  $C_c (G)$ i.e. a net $(\varphi_\alpha)$ in $C_c (G)$ with
$\varphi_\alpha \ast \psi \to \psi$ for all $\psi \in C_c (G)$. We
can additionally assume that for any relatively compact neighborhood
$U$ of $e\in G$ the supports of the $\varphi_\alpha$ are contained
in $U$ for sufficiently large $\alpha$. Now,  define $f_\alpha \in
L^2 (X,m)$ by
$$f_\alpha (\mu) := \nfunction_{\varphi_\alpha} (\mu) =  \int
\varphi_\alpha(s) d\mu (s).$$ Then, convergence
$$\nfunction^{f_\alpha} (\varphi) \to \nfunction_\varphi \mbox{ in $L^2 (X,m)$} $$
holds for all $\varphi \in C_c (G)$.  Indeed, by construction and
the defining properties of an approximate unit,   we have pointwise
(in $\mu \in X$) convergence
\begin{eqnarray*}
\nfunction^{f_\alpha}(\varphi) (\mu)& =& \int_G \varphi(t) \int_G
\varphi_\alpha (s - t) d\mu(s) dm_G (t)\\
&=& \int_G \int_G \varphi (t) \varphi_\alpha (s - t) dm_G (t) d\mu
(s)\\
&=& \int_G \varphi \ast \varphi_\alpha (s) d\mu (s)\\
&\to & \int_G \varphi (s) d\mu (s)\\
&=&\nfunction_\varphi (\mu)
\end{eqnarray*}
 as well as a uniform (in $\mu \in X$) bound
$$|\nfunction^{f_\alpha} (\varphi) (\mu)|\leq \|f\|_\infty |\mu| (K) \leq C
<\infty$$ due to the assumption on $X$. Thus, the maps
$\nfunction^{f_\alpha}$ converge to the map $\nfunction$. Moreover,
the spectral measures $\mu_{f_\alpha}$ converge to
$\widehat{\gamma}$ (see Corollary 1 in \cite{BL}).

\medskip

For our subsequent considerations a  special class of TMDS will be
particularly relevant. These are introduced next. Consider the set
$C(G)$ of continuous functions on $G$ equipped with the topology of
uniform local convergence. In particular, the map
$$\delta_0 : C(G)\longrightarrow \CC, f\mapsto f(e),$$
is continuous. Clearly, the group $G$ acts continuously on $C(G)$
via translations
$$G\times C(G)\longrightarrow C(G), (s,f)\mapsto f(s\cdot).$$
Whenever $Y$ is a compact and invariant subset of $C(G)$ and $m$ an
invariant probability measure on $Y$,  we call $(Y,G,m)$ a
\textit{function dynamical system} or just FDS for short. Let now an
FDS $(Y,G,m)$ be given. As $\delta_0$ is continuous, its range on
such a $Y$ is compact. As $Y$ is invariant,  we easily find  that
any $f\in Y$ must be bounded. Moreover, it is not hard to check that
the map
$$Y \longrightarrow M^\infty(G), f\mapsto f m_G,$$
is continuous and injective. In this way, any $(Y,G,m)$ can be
considered as a TMDS. Hence, FDS are special instances of TMDS.

\begin{remark} (a) For discrete groups $G$, the set $C(G)$ can
easily be identified with the set of all measures on $G$. Under this
identification any TMDS can then be seen as a FDS (and vice versa).
So, in this case there is no difference between FDS and TMDS.

(b) As discussed above   any FDS can be seen as a TMDS, whose
measures all have continuous densities with respect to Haar measure.
The converse is not true. More specifically, it is not hard to
construct examples of TMDS, all of whose elements have continuous
densities, which are not FDS. Consider e.g. on the real numbers
$\RR$ the measure $\mu = h \lambda$  with $\lambda$ the Lebesgue
measure and the density $h: \RR\longrightarrow [0,\infty)$ given by
$$ h:= \sum_{n\in\ZZ} |n| 1_{[n- 2^{-|n|}, n + 2^{-|n|}]},$$
where $1_A$ denotes the characteristic function of $A$. Then, the
measure  $\mu$ is a translation bounded measure. In fact, for any
$\delta >0$ it  even satisfies $\mu([s -\delta,s+ \delta]) \to 0,
|s|\to \infty$. This easily gives that the \textit{hull} of $\mu$
given by
$$\Omega (\mu):=\overline{\{t\mu : t\in \RR\}}$$
consists of the translates $t\mu$, $t\in\RR$, as well as the zero
measure. Then,  the hull is a TMDS and all of its elements have
densities with respect to Lebesgue measure. Moreover, the measures
$n\mu$, $n\in\NN$,  converge clearly to the zero measure. However,
the densities of $n\mu$ do not converge pointwise let alone locally
uniformly (but rather explode at $0$).
\end{remark}

\section{Spectral theory as diffraction theory of
factors}\label{sec-factors} In this section we will discuss how the
spectral theory of arbitrary dynamical systems can be considered as
diffraction theory of suitable factors.  The basic idea is the
following:  Whenever $(X,G,m)$ is a dynamical system and $f$ is a
continuous function on $X$, any $x\in X$ will encode a bounded
continuous function
$$\Phi^f_x : G\longrightarrow \CC \mbox{ via } \Phi_x^f (t) :=f(\alpha_{-t} x).$$
and $\Phi^f$ can be seen as a factor map (see below for details and
definitions). This effectively allows us to consider the situation
of $(X,G,m)$ together with an $f\in C(X)$ as an FDS. Now, as
discussed in the previous section, for such systems (and even more
general TMDS) diffraction theory has been developed in \cite{BL} and
the theory developed there can be applied to

\begin{itemize}

\item identify the autocorrelation $\gamma^f$ as the autocorrelation
of a function dynamical system, which is a  factor of the original
system;

\item exhibit the spectral measure $\varrho^f$ as the diffraction of
a a factor of the original system;

\item show pointwise existence of the autocorrelation along sampling
of an orbit.
\end{itemize}

Consider a dynamical system $(X,G,m)$. A  dynamical system $(Y,G,n)$
is called a \textit{factor} of $(X,G,m)$ if there exists a
continuous surjective $G$-invariant  map $\Phi : X\longrightarrow Y$
with $\Phi (m) = n$. Here, $\Phi (m)$ is the push-forward of the
measure $m$ by $\Phi$ defined by
$$\int_{Y} F(y) d\Phi(m) (y) = \int_X (F\circ \Phi) (x)\,  dm(x)$$
for measurable $F$ on $X$ with $F\geq 0$. In this situation the map
$\Phi$ is called \textit{factor map}.

Our first aim is to identify all the FDS factors of $(X,G,m)$. To do
so, we first consider a continuous  $f$ on $X$. Then,  we can indeed
define for any $x\in X$ the function $\Phi^f_x$  on $G$ via
$\Phi^f_x (t) = f(\alpha_{-t} x)$. With the measure $\mu^f_x :=
\Phi^f_x dt$ we then obviously have
$$\int \varphi (t) d\mu^f_x (t) = \int \varphi (t)   f(\alpha_{-t} x)
dm_G (t) = \nfunction^f (\varphi) (x)$$ for $\varphi \in C_c (G)$.
It is then not hard to see that
$$\Phi^f : X\longrightarrow C(G), x\mapsto \Phi^f_x,$$
is continuous and equivariant. Denote the image of $X$ under
$\Phi^f$ by $X^f$ and let $\Phi^f (m)$ be the  push-forward of $m$.
Then, $(X^f, \alpha, \Phi^f (m))$ is a FDS. It can easily be seen to
be ergodic if $(X,G,m)$ is ergodic. Conversely, whenever $\Phi :
X\longrightarrow Y$ is the factor map between $(X,G,m)$ and the FDS
$(Y,G,n)$, the function $f = \delta_0 \circ \Phi$ is continuous with
$\Phi = \Phi^f$.

We summarize the preceding considerations in the following
proposition.

\begin{prop}[Correspondence between FDS factors and continuous
$f$]\label{prop-correspondence} Let $(X,G,m)$ be a dynamical system.
For any continuous $f: X\longrightarrow \CC$ the FDS  $(X^f, G,
\Phi^f (m))$ is a factor of $(X,G,m)$ via the factor map $\Phi^f$
and the map
$$C (X)\longrightarrow \mbox{FDS factors of $(X,G,m)$}, \;  f\mapsto
(X^f,G,\Phi^f (m)),$$ is a bijection (with inverse given by $\Phi
\mapsto \delta_0\circ \Phi$).
\end{prop}

Having identified the FDS factors we now turn to study their
autocorrelation and diffraction. As is not surprising,  the
autocorrelation of $(X^f, G,\Phi^f (m))$ and of $f$ are identical.

\begin{lemma}[Autocorrelation of $f$ as autocorrelation of a
FDS]\label{lem-auto-TMDS} Let $(X,G,m)$ be a dynamical system and
$f: X\longrightarrow \CC$ continuous.  The autocorrelation
$\gamma^{(\Phi^f (m))}$  of $(X^f, \alpha, \Phi^f (m))$ is just
$\gamma^f$.
\end{lemma}
\begin{proof}
From the definition of $\gamma^f$ and of $\gamma^{(\Phi^f (m))}$ we
find for any $\varphi,\psi \in C_c (G)$
\begin{eqnarray*}
\gamma^{(\Phi^f (m))}(\varphi \ast \tilde{\psi})&=& \langle
\nfunction_\varphi, \nfunction_\psi\rangle_{L^2 (X^f, \Phi^f
(m))}\\&=& \int_{X^f} y (\varphi) \overline{y(\psi)} d\Phi^f (m) (y)\\
&=& \int_X \Phi^f_x (\varphi) \overline{\Phi^f_x (\psi)}  dm (x)\\
&=& \langle \nfunction^f (\varphi), \nfunction^f (\psi)\rangle_{L^2
(X,m)}\\
&=&\gamma^f (\varphi \ast \tilde{\psi})
\end{eqnarray*}
This finishes the proof.
\end{proof}

The lemma effectively allows one to reduce the spectral theory (of
continuous functions) of dynamical systems to diffraction theory for
TMDS. In particular, it is an main step in the proof of the
following result.

\begin{theorem}[Spectral theory as diffraction theory of FDS factors and vice versa]
\label{t:spectralmeasure-as-diffraction-of-TMDS} Let $(X,G,m)$ be a
dynamical system. Then the following holds.
\begin{itemize}
\item[(a)] For any FDS factor $(Y,G,n)$ with factor map $\Phi$ the equality
$\widehat{\gamma^{n}}  = \varrho^{\delta_0 \circ \Phi}$ holds.

\item[(b)] For any continuous $f : X\longrightarrow \CC$ the equality
$\widehat{\gamma^{\Phi^f (m)}}  = \varrho^f$ holds.

\end{itemize}
In particular, the map
\begin{eqnarray*}
\mbox{FDS factors of $(X,G,m)$} &\longrightarrow &\mbox{Spectral
measures of continuous functions}\\
 (Y,G,n) &\mapsto&
\widehat{\gamma}^{n}, \end{eqnarray*}
is onto.

\end{theorem}
\begin{proof} The  preceding lemma and Theorem
\ref{thm-main-abstract} directly give (b).  Now,  (a) is a direct
consequence of (b) and Proposition \ref{prop-correspondence}.  The
last statement is then  immediate from (a) and (b).
\end{proof}

\begin{remark} The theory of TMDS was originally developed in \cite{BL} with a view
towards dynamical systems of Delone sets and generalizations thereof
(as these are the relevant models for aperiodic order). The previous
considerations and, in particular, the previous Theorem show that
TMDS (and in fact even FDS) arise very naturally in the treatment of
spectral theory of \textit{any} dynamical system. Somewhat
informally, the previous theorem can be summarized as giving that
spectral theory of continuous functions is the same as diffraction
theory of FDS factors. On the structural level this can be seen as a
main result of the article.  \exend
\end{remark}

We are now heading towards a version of the previous result taking
simultaneously all spectral measures into account. To give this
version we need some further concepts.

Consider a dynamical system $(X,G,m)$ with associated unitary
representation $T$ of $G$. By Stone's theorem, compare
\cite[Sec.~36D]{Loomis}, there exists a (unique) projection-valued
measure
\begin{equation*}
  E_T\! : \; \{\mbox{Borel sets of $\widehat{G}$}\} \; \longrightarrow \;
  \{\mbox{projections on $L^2 (X,m)$}\}
\end{equation*}
with $$\langle f \ts | \ts E_T( \cdot) f\rangle = \varrho^f$$ for
any $f\in L^2 (X,m)$. Then, the projection-valued measure $E_T$
contains the entire spectral information on the dynamical system. It
is desirable to encode this spectral information in terms of
measures on $\widehat{G}$. One way of doing so is via the family of
spectral measures. More generally, we introduce the following
definition.

\begin{definition}[Complete spectral invariant]\label{spectralinvariant}
Let $T$ be the unitary representation associated to  $(X,G,m)$, and
let $E_T$ be the corresponding projection-valued   measure.  A
family $ \{ \sigma_\iota \}$ of measures on $\widehat{G}$ (with
$\iota$ in some index set $J$) is called a \emph{complete spectral
invariant} when $E_T (A) = 0$ holds for a Borel set
  $A\subset \widehat{G}$ if and only if $\sigma_\iota (A) = 0$ holds
  for all $\iota \in J$.
\end{definition}

\begin{theorem}[Diffractions of FDS factors as complete spectral invariant]\label{t:analogue-BLE} Let $(X,G)$ be an arbitrary dynamical system with
invariant probability measure $m$. Then, the diffractions of its FDS
factors form a  complete spectral invariant.
\end{theorem}
\begin{proof}
This is a direct consequence of Theorem
\ref{t:spectralmeasure-as-diffraction-of-TMDS} and the denseness of
$C(X)$ in $L^2 (X,m)$.
\end{proof}

\begin{remark} The main abstract result of \cite{BLE} deals with special TMDS viz
Delone dynamical systems with finite local complexity. It show that
for such systems  the diffraction of factors, which are also Delone
dynamical systems of finite local complexity, provide a complete
spectral invariant. The methods of \cite{BLE} heavily rely on the
assumption of finite local complexity. It it not even clear how to
extend them to general TMDS. The previous result now provides an
analogue  for arbitrary topological dynamical systems.
\end{remark}

\bigskip

So far, we have been concerned with FDS factors. It turns out that
it is also possible to work with  general TMDS factors.  Here, the
crucial insight is  that  the diffraction of TMDS factors is also
encoded by spectral theory. Here is the corresponding result.

\begin{prop}[Diffraction of TMDS factors via spectral measures]
\label{prop:diffraction-factor} Let $(X,G,m)$ be dynamical system
and  $(Y,G,\Phi (m))$ be a TMDS factor with factor map $\Phi :
X\longrightarrow Y$. Let $\widehat{\gamma}$ be the diffraction of
$(Y,G,\Phi(m))$. Then,
$$|\check{\varphi}|^2 \widehat{\gamma} = \varrho^{\nfunction_\varphi \circ
\Phi}$$ for all $\varphi \in C_c (G)$.
\end{prop}
\begin{proof} In this proof we will deal with both spectral measures
of functions on $X$ and on $Y$. In order to avoid confusion we
therefore write $\varrho^f_X$ to denote the  spectral measure of a
function $f\in C(X)$ and $\varrho^g_Y$ to denote the spectral
measure of a function $g\in C (Y)$. Similarly, we will write $T_t^X$
and $T_t^Y$ to denote the unitary action of $G$ on $L^2(X,m)$ and
$L^2 (Y,\Phi (m))$ respectively.

As discussed above we have
$$|\check{\varphi}|^2 \widehat{\gamma} =
\varrho^{\nfunction_\varphi}_Y$$ for all $\varphi \in C_c (G)$. As
$\Phi$ is a factor map, the definition of $\Phi (m)$ easily gives
that
$$\langle \nfunction_\varphi \circ \Phi, T_t^X \nfunction_\varphi
\circ \Phi\rangle_{L^2 (X,m)} = \langle \nfunction_\varphi, T_t^Y
\nfunction_\varphi\rangle_{L^2 (Y,\Phi (m))}$$ for all $t\in G$. By
the defining properties of the spectral measures, this yields
$$\varrho^{\nfunction_\varphi \circ \Phi}_X  = \varrho^{\nfunction_\varphi}_Y$$
for all $\varphi \in C_c (G)$. Putting this together we find
$$|\check{\varphi}|^2 \widehat{\gamma} =
\varrho^{\nfunction_\varphi}_Y =\varrho^{\nfunction_\varphi \circ
\Phi}_X$$ and the proof is finished.
\end{proof}

From the preceding considerations we now obtain the following
corollary.

\begin{coro}[Diffraction of TMDS factors as complete spectral invariant]\label{cor:analogue-BLE} Let $(X,G)$ be an arbitrary dynamical system with
invariant probability measure $m$. Then, the diffractions  of its
TMDS factors form  a complete spectral invariant.
\end{coro}
\begin{proof} This is a direct consequence of Theorem
\ref{t:analogue-BLE}  and Proposition \ref{prop:diffraction-factor}.
\end{proof}

As a consequence of Lemma \ref{lem-auto-TMDS}  it is also  possible
to  show that the autocorrelation exists as a pointwise limit along
the orbits of a sampling function. To do so, we need some further
concepts. Firstly, recall from \cite{Martin2} that a sequence
$(B_n)$ of compact subsets of $G$ is called a \textit{van Hove
sequence} if
\begin{equation*}
   \lim_{n\to \infty} \frac{|\partial^K B_n|}{|B_n|} \; = \; 0
\end{equation*}
for all compact $K\subset G$. Here, for compact $B,K$, the
``$K$-boundary''  $\partial^K B $  of $B$ is defined as
\begin{equation*}
   \partial^K B \; := \; \overline{((B + K)\setminus B)} \cup
    ((\overline{G\setminus B} - K) \cap B),
\end{equation*}
where the bar denotes the closure. The existence of van Hove
sequences for all $\sigma$-compact LCA groups is shown in
\cite[p.~249]{Martin2}, see also Section 3.3 and Theorem (3.L) of
\cite[Appendix]{Tempelman}. Moreover, every van Hove sequence is a
F{\o}lner sequence, i.e., $\lvert B_n \triangle (B_n + K) \rvert /
\lvert B_n \rvert \xrightarrow{n\to\infty} 0$, for every compact set
$K\subset G$, where $A\triangle B := (A \setminus B) \cap
(B\setminus A)$.

\begin{theorem}[Sampling for continuous functions]\label{sampling-cont} Let $(X,\alpha,m)$ be ergodic and
assume that $G$ has a countable basis of topology. Let  $(B_n)$
   be a van Hove sequence along which the Birkhoff ergodic theorem
   holds. Let $f : X\longrightarrow \CC$ be continuous.
Then, the measures $$\gamma_{x,n}^f :=\frac{1}{m_G (B_n)}
\mu^f_x|_{B_n} \ast \widetilde{\mu^f_x|_{B_n}}$$ converge vaguely to
$\gamma^f$ for almost every $x\in X$.
\end{theorem}
\begin{proof}
From  \cite{BL} we obtain that the limit exists almost everywhere
and equals the autocorrelation of the TMDS $(X^f, \alpha, \Phi^f
(m))$. By the Lemma \ref{lem-auto-TMDS}  that autocorrelation is
just $\gamma^f$ and the theorem is proven.
\end{proof}

\section{Application to  discrete spectrum}\label{s:application}
In this section we will use the material of the previous section to
study  discrete spectrum. In particular, we will combine this
material with a recent result of \cite{SpStr} to  provide a
characterization of discrete spectrum. This characterization is a
variant of well-known  characterizations for action of $\ZZ$ (see
end of the section for details).

\medskip

Recall that the dynamical system $(X,G,m)$ is said to have \textit{
discrete spectrum} if all spectral measures are pure point measures.
 We start with the following consequence of our discussion.
\begin{lemma}\label{l:char-pp-abstract} Let $(X,G,m)$ be a dynamical system. Then the
following assertions are equivalent.
\begin{itemize}
\item[(i)] The dynamical system $(X,G,m)$ has  discrete
spectrum.
\item[(ii)] For every $f\in C(X)$ the dynamical system $(Y_f,
G, \Phi^f (m))$ has  discrete spectrum.
\end{itemize}
\end{lemma}
\begin{proof} By denseness of $C(X)$ in $L^2 (X,m)$ discrete
spectrum is equivalent to $\varrho^f$ being a pure point measure for
all $f\in C(X)$. By Theorem
\ref{t:spectralmeasure-as-diffraction-of-TMDS} this is equivalent to
pure pointedness of all measures $\widehat{\gamma}^{\Phi^f (m)}$ for
$f\in C(X)$. Now, for TMDS  a main result of \cite{BL} gives that
pure pointedness of the diffraction spectrum is equivalent to
discrete spectrum of the dynamical system. This yields the desired
statement.
\end{proof}

\begin{remark} Other features of the dynamical system (mean equicontinuity,
equicontinuity etc) can be characterized similarly via factors. We
refrain from giving details.
\end{remark}

We need some further concepts. Consider an arbitrary van-Hove
sequence $(B_n)$. Then, we define the associated mean via
$$ \overline{M}:=\overline{M}_{\mathcal{B}}  : \mbox{Bounded function on G}
\longrightarrow [0,\infty), \overline{M} (h) :=\limsup_{n\to \infty}
\frac{1}{m_G(B_n)} \int_{B_n} h(s) d m_G (s).$$ Clearly,
$\overline{M}$ is a seminorm on the space of bounded functions on
$G$. A subset $A$ of $G$ is called \textit{relatively dense} if
there exists a compact set $K\subset G$ with
$$ G =\bigcup_{a\in A} (a + K).$$
Now, a uniformly continuous bounded $h : G\longrightarrow \CC$ is
called \textit{mean almost periodic} if for every $\varepsilon
>0$ the set
$$\{t\in G : \overline{M}_{\mathcal{B}} ( h - h(t\cdot))
<\varepsilon\}$$ is relatively dense in $G$. A translation bounded
measure $\mu$ on $G$ is called \textit{mean almost periodic} if for
every $\varphi \in C_c (G)$ the function $\mu\ast \varphi$ is mean
almost periodic. Then, it is not hard to see that a uniformly
continuous bounded function $h : G\longrightarrow \CC$ is mean
almost periodic if and only if the associated measure $h m_G$ is
mean almost periodic (compare Proposition 3.9  in \cite{SpStr}).

\begin{theorem} Let $(X,G,m)$ be an ergodic dynamical system. Assume that $G$ has a countable basis of the topology and
 let $(B_n)$ be a van Hove sequence along which the
Birkhoff ergodic theorem holds. Then, the following assertions are
equivalent.
\begin{itemize}
\item[(i)] The dynamical system $(X,G,m)$ has discrete
spectrum.
\item[(ii)] For every $f\in C(X)$
the function $\Phi^f_x$ is mean almost periodic for $m$ almost every
$x\in X$.
\end{itemize}
\end{theorem}
\begin{proof} By the preceding lemma,  (i) is equivalent to
$(Y_f, G, \Phi^f (m))$ having  discrete spectrum for all $f\in
C(X)$.  Now, $(Y_f, G, \Phi^f (m))$ is a TMDS (and even an FDS) over
a locally compact abelian group with countable basis of the
topology.  Thus, by a main result of \cite{SpStr} discrete spectrum
of $(Y_f, G, \Phi^f (m))$ is equivalent to mean almost periodicity
of $y m_G$ for  $\Phi^f (m)$ almost every  $y$. As discussed above
this is equivalent to mean almost periodicity of the function $y$
for $\Phi^f (m)$ almost every $y\in Y_f$.  Now, the desired
statement follows from the definition of $\Phi^f(m)$.
\end{proof}

Let us conclude this section with a consequence of the preceding
theorem in the case that  $X$ has a countable basis of the topology.
In this case  the set $C(X)$ is separable in the uniform topology.
Thus, one can then  chose  in (ii) one set $X'$ of full measure in
$X$ such that $\Phi^f_x$ is mean almost periodic for every $x\in X'$
and all $f\in C(X)$.  If one now calls  a point  $x\in X$
\textit{mean almost periodic} if $\Phi^f_x$ is mean almost periodic
for every $f\in C(X)$, the previous result can then be restated as
follows:

\begin{coro} Let $(X,G,m)$ be an ergodic dynamical system.
Assume that both $X$ and $G$ have a countable basis of the topology.
Then, the dynamical system $(X,G,m)$ has discrete spectrum if and
only if almost every $x\in X$ is mean almost periodic.
\end{coro}

\begin{remark} We could not find the result in this form in the
literature. However, it will  certainly not be a surprise to
experts.  Indeed, for ergodic systems over  $G = \ZZ$, it is known
that discrete spectrum is equivalent to $\ZZ\ni n\mapsto f(nx)$
belonging to the Besicovich class for almost every $x\in X$ whenever
$f$ is a bounded measurable function on $X$, see  Theorem 3.22 in
\cite{BeLo}. The condition of Besicovich class seems to be somewhat
stronger than mean almost periodicity. On the other hand, clearly,
the condition that $f$ is continuous is stronger than $f$ being
bounded and measurable. We also note that a related result for
subshifts over a finite alphabet  can be found in Lemma 5 of 
\cite{Ver}. There, discrete spectrum is characterized  with a mean
almost periodicity condition on points rather than functions.
\end{remark}

\section{The measurable situation}\label{sec-measurable}
The previous sections has been set in the topological category.
Starting with a dynamical system $(X,G,m)$ we considered continuous
functions on $X$ as well as  factor maps. From the point of view of
the unitary representation $T$ it is also desirable to treat
functions $f\in L^2 (X,m)$. Roughly speaking this case can be
treated by similar means as the case of continuous $f$ by using the
theory of square integrable measure dynamical systems from \cite{LS}
instead of the theory of translation bounded dynamical systems from
\cite{BL}. A sketch is given in this section.

\bigskip

Recall that we denote by $M (G)$ the set of all measures on $G$.
This set is equipped with the vague topology and the associated
Borel-$\sigma$-algebra. There is a canonical action of $G$ by
translation. Any $\varphi \in C_c (G)$ then gives rise to a
continuous function $\nfunction_\varphi$ on $M(G)$ via
$$\nfunction_\varphi (\mu) := \int_G \varphi \, d\mu.$$
An invariant probability measure $n$ on $M(G)$ is called
\textit{square integrable} if
$$\int_{M(G)} |\nfunction_\varphi (|\mu|)|^2 dn (m) <\infty$$
for all $\varphi \in C_c (G)$. As shown in \cite{LS} any square
integrable measure $n$  then comes with a unique measure $\gamma =
\gamma^n $ satisfying
$$\gamma (\varphi \ast \widetilde{\psi}) = \langle \nfunction_\varphi, \nfunction_\psi
\rangle_{L^2 (M(G),n)}$$ for all $\varphi,\psi \in C_c (G)$. This
measure $\gamma$ is  called the \textit{autocorrelation of $n$}. Its
Fourier transform exists and is a positive measure called the
\textit{diffraction of $n$}.  Note that this generalizes the
framework of TMDS discussed above.

We now want to proceed as above by defining the measure $\mu^f_x$ on
$G$ via
$$\mu^f_x (\varphi) = \int_G \varphi(s) f(\alpha_{-s} x) d m_G
(s).$$ Here, we have to overcome the obstacle that this will not
necessarily make sense (due to possible unboundedness of $f$).
 However, it
turns out that it will make sense for almost every $x\in X$ under
suitable additional assumptions.

\begin{prop} \label{square-int} Let $f\in L^2 (X,m)$ be given. Then,
$$ \int_X \left(\int \varphi (t) |f(\alpha_{-t} x)| dm_G (t)\right)^2
dm (x) <\infty$$ for any $\varphi \in C_c (G)$ with $\varphi \geq
0$. In particular,  for $m$-almost every $x\in X$ there exists a
measure  $\mu^f_x$   $G$ with
$$\int \varphi (t) d\mu^f_x (t) =
\int \varphi (t)   f(\alpha_{-t} x) dm_G (t)$$ for $\varphi \in C_c
(G)$
\end{prop}
\begin{proof}
The first inequality  follows easily from Cauchy--Schwarz inequality
and Fubini's Theorem.  As for the last statement we note that
 by $\sigma$-compactness of
$G$ there exists a increasing  sequence of compact sets $K_n$ with
$$\bigcup_n \mbox{int}(K_n) = G.$$
Here, $\mbox{int}$ denotes the interior of a set.   For each $K_n$
we can now chose a function $\varphi_n \in C_c (G)$  with $0 \leq
\varphi_n \leq 1$ and $\varphi_n = 1$ on $K_n$. By the first
statement, we can then find a set $N_n\subset X$ with $m(N_n) =0$
such that
$$\int \varphi_n (t) |f(\alpha_{-t} x)| dm_G
(t) <\infty$$ for all $x\notin N_n$. Hence, we infer that $\mu^f_x$
is a measure for all $x\notin \bigcup N_n$ and this is the desired
statement.
\end{proof}

We denote the set of $x\in X$ for which $\mu^f_x$ is a measure by
$X_{meas}$. By the previous proposition this set has full measure.
Moreover, it is clearly invariant.   Then the map
$$\mu^f : X_{meas} \longrightarrow M (G)$$
is measurable. Denote by $\mu^f (m)$ the push-forward of the measure
$m$ to $M(G)$. This is an invariant probability  measure on $M(G)$
and ergodic if $m$ is ergodic. From Proposition \ref{square-int} we
easily infer that $\mu^f (m)$ is square integrable. Hence, we can
apply the theory of \cite{LS} to $\mu^f (m)$ to conclude existence
of $\gamma^{\mu^f(m)}$. Then, arguing as above for the case of TMDS
we find
$$\gamma^{\mu^f  (m)} = \gamma^f \mbox{ and } \widehat{ \gamma^{\mu^f  (m)}} =
\varrho^f.$$ Moreover, using the results of \cite{LS} we can then
also conclude (as above) that Theorem \ref{sampling-cont} remains
valid for $f\in L^2 (X,m)$.

\begin{remark} It may be worthwhile to emphasize the following
application of the results of this section: Any measure appearing as
spectral measure of a dynamical system can  be obtained as
diffraction measure of a whole family of functions. In principle,
this opens up a way to generate rather naturally diffraction
measures with interesting properties.
\end{remark}

\begin{appendix}

\section{A short review of diffraction theory}
In this section we briefly review the framework for diffraction
theory developed in \cite{LM}.

\bigskip

The basic pieces of data are given by a dynamical system $(X,G,m)$
(with $G$ a locally compact, $\sigma$-compact abelian group) and a
map
$$\nfunction : C_c (G) \longrightarrow L^2 (X,m)$$
with the following properties\footnote{The article \cite{LM}
introduced one more assumption, viz that $\nfunction$ is real. This,
however, is not needed for the results we discuss here.}:

\begin{itemize}

\item[(N1)] $\nfunction $ is linear.


\item[(N2)] $\nfunction$ is equivariant i.e. satisfies $\nfunction  (\varphi
(\cdot - t)) = \nfunction (\varphi) (\alpha_{-t}(\cdot))$ for all
$\varphi \in C_c (G)$ and all $t\in G$.

\item[(N3)] There exists a measure $\gamma$ on $G$ with  $\gamma (\varphi \ast \widetilde{\psi}) = \langle \nfunction (\varphi), \nfunction
(\psi)\rangle_{L^2 (X,m)}$ for all $\varphi,\psi \in C_c (G)$.

\end{itemize}

In this case $\gamma$ is a positive definite measure called the
\textit{autocorrelation} of $\nfunction$ and its Fourier transform
is a positive measure called the \textit{diffraction} of
$\nfunction$. Moreover, there exists a unique isometry, called
\textit{diffraction to dynamics map}
$$\Theta : L^2 (\widehat{G},\mu_f)\longrightarrow L^2 (X,m)$$
with $\Theta^f (\widehat{\sigma}) = \nfunction^f (\sigma)$ for any
$\sigma \in C_c (G)$.

It turns out that in applications one can often conclude validity of
(N3)  for maps satisfying (N1) and  (N2)  from  a continuity
property.
  Specifically, the
 following is shown in \cite{LM}.

 \begin{lemma} Let $(X,G,m)$
 be a dynamical system. Let $\nfunction : C_c (G) \longrightarrow L^2 (X,m)$
 satisfy {\textrm{(N1), (N2)}} and $\overline{\nfunction
(\varphi)} = \nfunction( \overline{\varphi})$ for all $\varphi \in
C_c (G)$.  If for all  $1\leq p,q\leq \infty$ with $1/p + 1/q = 1$
be given and any compact $K\subset G$ there exists a  $C_K>0$ with
$$|\langle \nfunction^f(\varphi),  \nfunction^f (\psi)\rangle| \leq C_K \|\varphi \|_{L^p (G) } \|\psi\|_{L^q (G)}$$
for all $\varphi,\psi\in C_c (G)$ with support contained in $K$ then
$\nfunction$ satisfies {\textrm{ (N3)}} as well.
\end{lemma}

\end{appendix}

\end{document}